\documentclass{article}
\usepackage{color}
\usepackage[colorlinks]{hyperref}
\usepackage{amsmath}
\usepackage{amssymb, amsthm}
\newtheorem{theorem}{Theorem}[section]
\newtheorem{remark}{Remark}

\newcommand{\N}{\mathbb{N}}
\newcommand{\ds}{\displaystyle}
\newcommand{\pf}{P$_{\textup V}$}
\newcommand{\pt}{P$_{\textup{III}}$}
\begin{document}
\title{
Recurrence coefficients of generalized Charlier polynomials and
the fifth Painlev\'e equation}
\author{Galina Filipuk, Walter Van Assche}
\date{\today}
\maketitle

\begin{abstract}
We investigate generalizations of the Charlier polynomials on the
lattice $\mathbb{N}$, on the shifted lattice $\mathbb{N}+1-\beta$
and on the bi-lattice $\mathbb{N}\cup (\mathbb{N}+1-\beta)$. We
show that the coefficients of the three-term recurrence relation
for the orthogonal polynomials are related to solutions of the
fifth Painlev\'e equation \pf\ (which can be transformed to the
third Painlev\'e equation). Initial conditions for different
lattices can be transformed to the classical solutions of  
\pf\ with special values of the parameters.
\end{abstract}



\section{Introduction}

In this paper we investigate the recurrence coefficients of
orthogonal polynomials, in particular generalized Charlier
polynomials, and show that they are related to the solutions of a
continuous Painlev\'e equation when viewed as functions of one of
the parameters. We are interested to see how the properties of the
orthogonal polynomials are related to properties of
transformations of the Painlev\'e equation.

The paper is organized as follows. In the introduction we shall
first review orthogonal polynomials  for a generalized Charlier
weight on different lattices and their main properties, following
\cite{bilattice}. Next we shall briefly recall the fifth
Painlev\'e equation and its B\"acklund transformation.  We then
derive the Toda system and show that the recurrence coefficients
can be expressed in terms of the solutions of  the fifth
Painlev\'e equation. Finally we study the initial conditions of
the recurrence coefficients for different lattices.

\subsection{Orthogonal polynomials for the generalized Charlier weight}

One of the most important properties of orthogonal polynomials is
a three-term recurrence relation. For a sequence $(p_n(x))_{n \in
\N}$ of orthonormal polynomials for a weight $w$ on the lattice
$\mathbb{N}=\{0,1,2,3,\ldots\}$
\begin{equation}\label{orthonormality}
\sum_{k=0}^{\infty} p_n(k)p_k(k)w(k)=\delta_{n,k},
\end{equation}
where $\delta_{n,k}$ is the Kronecker delta,
 this relation takes the following form:
\begin{equation}\label{3 term orthonormal}
 x p_n(x)=a_{n+1}p_{n+1}(x)+b_np_n(x)+a_np_{n-1}(x).
\end{equation}
 For the   monic
polynomials $P_n$ related to orthonormal polynomials
$p_n(x)=\gamma_n x^n+\cdots $ with
$$\frac{1}{\gamma_n^2}=\sum_{k=0}^{\infty}P_n^2(k)w(k), $$
the recurrence relation is given by
$$xP_n(x)=P_{n+1}(x)+b_nP_n(x)+a_n^2 P_{n-1}(x).$$

The (classical) Charlier polynomials (\cite[Chapter VI]{Chihara},
\cite{bilattice}) are orthogonal on the lattice $\mathbb{N}$ with
respect to the
Poisson distribution: %
$$%
 \sum_{k=0}^{\infty}C_n(k;a)C_m(k;a)\frac{a^k}{k!}=a^{-n}\exp(a)\, n!\,\delta_{n,m},\qquad a>0.%
$$ The weight $w_k=w(k)=a^k/k!$ satisfies the Pearson equation
$$\nabla w(x)=\left(1-\frac{x}{a}\right)w(x) ,$$ where $\nabla$ is the backward difference operator 
$$\nabla f(x)=f(x)-f(x-1).$$ 
Here the
 function $w(x)=a^x/\Gamma(x+1)$ gives the weights $w_k=w(k)=a^k/k!$. The
 Pearson equation for the Charlier polynomials is, hence, of the form
\begin{equation}\label{Pearson}\nabla[\sigma(x)w(x)]=\tau(x)w(x)
\end{equation}
with $\sigma=1$ and $\tau$ a polynomial of degree 1.

  Recall that the classical
 orthogonal polynomials are characterized by the Pearson equation
 (\ref{Pearson})
 with $\sigma$ a polynomial of degree at most 2 and $\tau$
 a polynomial of degree 1. Note that in (\ref{Pearson}) the
 operator
 $\nabla$ is used for orthogonal polynomials on the lattice and
   it is replaced by differentiation in case of orthogonal
  polynomials on an interval of the real line. The Pearson equation plays an
 important role in the theory of classical orthogonal polynomials since it
 allows to find many useful properties of these polynomials. 
It is known that   the recurrence coefficients of the Charlier
polynomials  are given explicitly  by
 $a_n^2=na$ and $b_n=n+a$ for $n\in\mathbb{N}.$

 The Charlier weight can be generalized \cite{bilattice}. One can
 use the weight function $$w(x)=\frac{\Gamma(\beta)
 a^x}{\Gamma(\beta+x)\Gamma(x+1)}$$ which gives the weight \begin{equation}\label{weight}
 w_k=w(k)=\frac{a^k}{(\beta)_k
 k!},\qquad a>0, \end{equation}
on the lattice $\mathbb{N}$. The
 Pearson equation is now of the form $$\nabla w(x)=\frac{a-x(\beta-1)-x^2}{a}\ w(x).$$
  The monic orthogonal polynomials $P_{n}(x;a,\beta)$ for the weight (\ref{weight})
 satisfy
\begin{equation}\label{orthogonality1}\sum_{k=0}^{\infty}P_n(k;a,\beta)P_m(k;a,\beta)\frac{a^k}{(\beta)_k
k! }=0,\qquad n\neq m,\end{equation} and the three-term recurrence
relation is
$$xP_{n}(x;a,\beta)=P_{n+1}(x;a,\beta) +b_n
P_{n}(x;a,\beta)+a_n^2 P_{n-1}(x;a,\beta), $$ with initial
conditions $P_{-1}=0 $ and $P_0=1.$ If we replace $a$ by $\beta a$
and then take the limit $\beta\rightarrow \infty$, one recovers
the Charlier polynomials \cite{bilattice}.

\begin{theorem}{\rm \cite[Th. 2.1]{bilattice}}\label{thm:
discrete} The recurrence coefficients for the orthogonal
polynomials defined by (\ref{orthogonality1}) for the weight
(\ref{weight}) on the lattice $\N$  satisfy
\begin{eqnarray}\label{discrete1}
&&b_n+b_{n-1}-n+\beta=\frac{an}{a_n^2}
\\ \label{discrete2}
&&(a_{n+1}^2-a)(a_n^2-a)=a(b_n-n)(b_n-n+\beta-1)
\end{eqnarray}
with initial conditions \begin{equation}\label{initial}
a_0^2=0,\;\;b_0=\frac{\sqrt{a}
I_{\beta}(2\sqrt{a})}{I_{\beta-1}(2\sqrt{a})},\end{equation} where
$I_{\beta}$ is the modified Bessel function.
\end{theorem}

The system (\ref{discrete1})--(\ref{discrete2})   can be
identified as a limiting case of an asymmetric discrete Painlev\'e
equation. In this paper we show that it can be obtained from the
B\"acklund transformation of the fifth Painlev\'e equation.

 One can use the
 weight (\ref{weight}) on the shifted lattice $\mathbb{N}+1-\beta$ and one can
 also combine both lattices to obtain the bi-lattice
 $\mathbb{N}\cup(\mathbb{N}+1-\beta)$. The orthogonality
 measure for the bi-lattice is a linear combination of the
 measures on $\mathbb{N}$ and $\mathbb{N}+1-\beta$.
The weight $w$  in (\ref{weight}) on the  shifted lattice
$\mathbb{N}+1-\beta=\{1-\beta,2-\beta,3-\beta,\ldots\}$ is given
by
\begin{equation}\label{weight sh}v_k:=w(k+1-\beta)=\frac{\Gamma(\beta)a^{1-\beta}}{\Gamma(2-\beta)}
\frac{a^k}{k!(2-\beta)_k},\qquad k\in\mathbb{N}. \end{equation} We
see that up to a constant factor it is equal to the weight on the
original lattice $\N$ with different parameters.
 The corresponding
monic orthogonal polynomials $Q_n$ on the shifted lattice satisfy
$$\sum_{k=0}^{\infty}Q_n(k+1-\beta)Q_{m}(k+1-\beta)v_k=0,\qquad n\neq m, $$
and one has $$Q_n(x)=P_n(x+\beta-1;a,2-\beta).$$ In
\cite[Th.~2.2]{bilattice} it is shown that  the recurrence
coefficients in the three-term recurrence relation $$x
Q_n(x)=Q_{n+1}(x)+ {\hat{b}}_nQ_n(x)+\hat{a}_n^2Q_{n-1}(x) $$
satisfy the same system (\ref{discrete1})--(\ref{discrete2}) (with
hats) but with initial conditions
\begin{equation}\label{initial sh}\hat{a}_0^2=0,\;\;{\hat{b}}_0=\frac{\sqrt{a}
I_{-\beta}(2\sqrt{a})}{I_{1-\beta}(2\sqrt{a})}.\end{equation}

Using the orthogonality measure $\mu=c_1\mu_1+c_2\mu_2,$ where
$c_1,\;c_2>0,$ $\mu_1$ is the discrete measure on $\mathbb{N}$
with weights $w_k=w(k)$ and $\mu_2$ is the discrete measure on
$\mathbb{N}+1-\beta$ with weights $v_k=w(k+1-\beta)$, one can
study the monic orthogonal polynomials
$R_n(x)=R_n(x;a,\beta,\tau),$ where $\tau=c_2/c_1>0$, satisfying
the three-term recurrence relation
$$x R_n(x)=R_{n+1}(x)+\tilde{b}_n
R_n(x)+\tilde{a}_n^2R_{n-1}(x).$$ The orthogonality relation is
given by
$$c_1\sum_{k=0}^{\infty}R_n(k)R_m(k)w_k+c_2\sum_{k=0}^{\infty}R_n(k+1-\beta)R_m(k+1-\beta)v_k=0,\qquad m\neq n. $$
According to \cite[Th.~2.3]{bilattice} the recurrence coefficients
$\tilde{a}_n^2$ and $\tilde{b}_n$ satisfy  the system
(\ref{discrete1})--(\ref{discrete2}) (with tilde) but with initial
conditions
\begin{equation}\label{initial
bi}\tilde{a}_0^2=0,\;\;{\tilde{b}}_0=\sqrt{a}\frac{
I_{\beta}(2\sqrt{a})+\tau
I_{-\beta}(2\sqrt{a})}{I_{\beta-1}(2\sqrt{a})+\tau
I_{1-\beta}(2\sqrt{a})}.
\end{equation}

Thus, it is shown in \cite{bilattice} that the orthogonal
polynomials for the generalized Charlier weight
 on the
lattice $\mathbb{N}$, on the shifted lattice $\mathbb{N}+1-\beta$
and on the bi-lattice $\mathbb{N}\cup (\mathbb{N}+1-\beta)$ have
recurrence coefficients $a_n^2$ and $b_n$ which satisfy the same
nonlinear system of discrete (recurrence) equations  but the
initial conditions are different for each case.

\subsection{The fifth Painlev\'e equation and its B\"acklund transformations}

The Painlev\'e equations possess the so-called Painlev\'e
property: the only movable singularities of the solutions are
poles \cite{GrLSh}. They are often referred to as nonlinear
special functions and have numerous applications in mathematics
and mathematical physics.

The fifth Painlev\'e equation \pf\ is given by
\begin{equation}\label{P5}
y''  =
\left(\displaystyle\frac{1}{2y}+\displaystyle\frac{1}{y-1}\right)(y')^2
-
\displaystyle\frac{y'}{t}+\displaystyle\frac{(y-1)^2}{t^2}\left(A
y + \displaystyle\frac{B}{y}\right)+\displaystyle\frac{C
y}{t}+\,\,\displaystyle\frac{D y (y+1)}{y-1},
\end{equation}
where $y=y(t)$ and $A,\; B,\; C,\;D$ are arbitrary complex
parameters.

The fifth Painlev\'e equation \pf\ in case $C\neq 0,\;\;D=0$ can
be reduced to the third Painlev\'e equation \cite{GrLSh}. There
exists a B\"acklund transformation between solutions of the fifth
Painlev\'e  equation with $D=0$.

\begin{theorem} {\rm \cite{Adler}, \cite{Tsegelnik}}\label{thm: BTr} Let
$y=y(t)$ be a solution of \pf\ with parameters $A,\,B,\,C,\,D$
given by $B=-p^2/2\neq 0$, $C\neq 0$ and $D=0$. Then the function
$y_1=y_1(t)$ given by $$y_1=1-\frac{2C t y^2(y-1)}{2A
y^2(y-1)^2-(p(y-1)-t y')^2}$$ is a solution of \pf\ with
parameters $A_1=A,\;B_1=-(p-1)^2/2,\;C_1=C,\;D_1=0$. The inverse
transformation is given by $$y=1-\frac{2C t y_1^2(y_1-1)}{2A
y_1^2(y_1-1)^2-((1-p)(y_1-1)-ty_1')^2}.$$
\end{theorem}

In case $D\neq 0$ there exists another well-known transformation.
\begin{theorem} {\rm \cite[Th. 39.1]{GrLSh}} \label{thm: BTr2}
If $y=y(t)$ is the solution the fifth Painlev\'e equation
(\ref{P5}) with parameters $A, \,B,\,C,\,D$, then the
transformation
$$T_{\varepsilon_1,\varepsilon_2,\varepsilon_3}:\;y\rightarrow
y_1$$  gives another solution $y_1=y_1(t)$ with new values of the
parameters $A_1, \,B_1,\,C_1,\,D_1$, where
\begin{eqnarray*}&&
y_1=1-\frac{2d\,t\,y}{ty'-ay^2+(a-b+d\,t)y+b},\\
&&A_1=-\frac{1}{16D}(C+d(1-a-b))^2,\quad
B_1=\frac{1}{16D}(C-d(1-a-b))^2,\\ && C_1=d(b-a),\;\;D_1=D,
\\
&&a=\varepsilon_1\sqrt{2A},\;\;b=\varepsilon_2\sqrt{-2B},\;\;d=\varepsilon_3\sqrt{-2D},\quad
\varepsilon_j^2=1,\;\;j\in\{1,2,3\}.
\end{eqnarray*}
\end{theorem}

{\bf Example.} Let $y=y_n(t)$ be a solution of \pf\ with
parameters $A=-B=n^2/8,\;\;C=0,\;\;D=-8.$ Then the function
$$y_{n+1}=1-\frac{16 t \varepsilon y}{2t y'+n y^2+8t \varepsilon
y-n }$$ is  a solution of \pf\ with
$A=-B=(n+1)^2/8,\;\;C=0,\;\;D=-8,\;\;\varepsilon^2=1.$

\begin{remark}
{\rm For our purposes it is sufficient to use the standard
B\"acklund transformations of the Painlev\'e transcendents which
are given in NIST Digital Library of Mathematical Functions  (DLMF
project\footnote{http://dlmf.nist.gov/32.7}), see \cite[Chapter 32]{NIST}. This also
contains various algebraic aspects of the Painlev\'e equations. The interested
reader can easily  re-formulate the transformations within the
framework of Noumi-Yamada's birational representation of the
affine Weyl groups \cite{Noumi}.}
\end{remark}

\subsection{Relation to the third Painlev\'e equation}

The third Painlev\'e equation \pt\  is given by
\begin{equation}\label{P3}
u''=\frac{u'^2}{u}-\frac{u'}{z}+\frac{1}{z}(\tilde{\alpha}
u^2+\tilde{\beta})+\tilde{\gamma} u^3-\frac{\tilde{\delta}}{u},
\end{equation}
where $u=u(z)$ and $\tilde{\alpha},\,\tilde{\beta},\,\tilde{\gamma},\,\tilde{\delta}$ are
arbitrary constants. By a scaling  of the independent variable, the
parameters $\tilde{\gamma}$ and $\tilde{\delta}$ can be taken equal to $1$ and
$-1$ respectively. There is a transformation between solutions of
\pf\ and \pt.

\begin{theorem}\cite[Th. 34.1]{GrLSh}\label{P3P5}
Let $u=u(z,\tilde{\alpha},\tilde{\beta})$ ($\tilde{\gamma}=1$, $\tilde{\delta}=-1$) be a solution
of \pt\ such that $R=u'-\varepsilon u^2-(\tilde{\alpha} \varepsilon
-1)u/z+1\neq 0.$ Then the function
$$y(t)=1-\frac{2}{R(\sqrt{2t})},\;\;\varepsilon^2=1,$$ is a solution of \pf\ with
parameters
$$A=\frac{(\tilde{\beta}-\tilde{\alpha}\varepsilon+2)^2}{32},\,B=-\frac{(\tilde{\beta}+\tilde{\alpha}\varepsilon-2)^2}{32},\,C=-\varepsilon,\,D=0.$$
\end{theorem}

The inverse transformation is given by the following theorem.
\begin{theorem}\cite[Th. 34.3]{GrLSh}\label{P5P3}
Let $y(t)$ be a solution of \pf\ with parameters
$A,\,B,\,C^2=1,\,D=0.$ Then the function
$$u(z)=\frac{\sqrt{2t}y}{\Phi(t)},\;\;z^2=2t,\;\;\Phi(t)=t y'-\sqrt{2A} y^2+(\sqrt{2A}+\sqrt{-2B})y-\sqrt{-2B}\neq
0$$is a solution of \pt\ with parameters
$$\tilde{\alpha}=2C(\sqrt{2A}-\sqrt{-2B}-1),\;\tilde{\beta}=2(\sqrt{2A}+\sqrt{-2B}),\;\tilde{\gamma}=1,\;\tilde{\delta}=-1.$$

\end{theorem}
\section{Main results}

In this section we  relate the recurrence coefficients to the
solutions of the fifth and third  Painlev\'e equations and
identify the discrete system (\ref{discrete1})--(\ref{discrete2}).

\subsection{The Toda system}

It is well-known that for a
measure of the form $e^{tx} \, d\mu(x)$ the recurrence
coefficients satisfy
\begin{equation} \label{Toda0}
\begin{cases} \ds \left(a_n^2\right)':=\frac{d}{dt}\left(a_n^2\right)=a_n^2(b_{n}-b_{n-1}), \\
\ds b_n':=\frac{d}{dt}b_n=(a_{n+1}^2-a_n^2).
\end{cases}
\end{equation}
with initial conditions $a_n^2(0)$ and $b_n(0)$ given by the
recurrence coefficients of the orthogonal polynomials for the
measure $\mu$. The weight (\ref{weight}) is of this form when we
put $a = e^t$. Hence, the recurrence coefficients of the
generalized Charlier polynomials on the lattice $\N$, when viewed
as functions of the parameter $a$ with respect to weight
(\ref{weight}), satisfy
\begin{equation} \label{Toda}
\begin{cases} \ds \left(a_n^2\right)':=\frac{d}{da}\left(a_n^2\right)=\frac{a_n^2}{a}(b_{n}-b_{n-1}), \\
\ds b_n':=\frac{d}{da}b_n=\frac{1}{a}(a_{n+1}^2-a_n^2).
\end{cases}
\end{equation}
We get this system from (\ref{Toda0}) by applying the
transformation $a=\exp(t)$ (see \cite{LiesGFWalter}). The weight (\ref{weight sh}) is also
of this form when $a = e^t$, since
$$   w(k+1-\beta) = \frac{\Gamma(\beta) e^{t(k+1-\beta)}}{\Gamma(2-\beta) k! (2-\beta)_k}.  $$
Hence, the recurrence coefficients for the shifted lattice also
satisfy the Toda system (\ref{Toda}). Also for the bi-lattice we
have that $w(x) = e^{tx} \Gamma(\beta)/\Gamma(\beta+x)
\Gamma(x+1)$ when $a=e^t$, which is of the required form, so that
the recurrence coefficients again satisfy the Toda system
(\ref{Toda}). The only difference between these three cases is in
the initial conditions which correspond to the recurrence
coefficients of the generalized Charlier polynomials with $n=0$  
on respectively the lattice $\N$, the shifted lattice $\N+1-\beta$
and the bi-lattice $\N \cup (\N+1-\beta)$.

\subsection{Relation to the B\"acklund transformation of the fifth  Painlev\'e equation}

First, let us obtain a nonlinear discrete equation for $b_n(a)$.
Denoting  $a_n^2=x_n$, we have from equation (\ref{discrete2})
  that
\begin{equation}\label{xn+1}
x_{n+1}=\frac{a(a+n(\beta-n-1)+(1+2n-\beta-b_n)b_n-x_n)}{a-x_n}.
\end{equation}
Equation (\ref{discrete1})   gives $$x_n=\frac{a
n}{\beta-n+b_{n-1}+b_n}.$$ Substituting the last equation and
(\ref{xn+1}) into equation (\ref{discrete1}) with $n+1$ we have a
nonlinear difference equation of the form
\begin{equation}\label{nonlinear discrete}
F(b_{n-1},b_n,b_{n+1},a)=0,
\end{equation}
for some $F$ which we do not write out explicitly, since it is
quite long. We shall show later on that equation (\ref{nonlinear
discrete}) can in fact be obtained from the B\"acklund
transformation of the fifth Painlev\'e equation.

To derive the differential equation for $b_n(a)$ we use
(\ref{xn+1}) and from equation (\ref{discrete1})   we have
$$b_{n-1}=\frac{a n+n x_n-\beta x_n-b_n x_n}{x_n}.$$ Substituting
this into the first equation of (\ref{Toda})   we get
$$x_n'=\frac{(\beta-n+2b_n)x_n-a n}{a}$$ and from the second
equation of (\ref{Toda}) we have
\begin{equation}\label{bn'}
b_n'=\frac{a(n-a+n^2-n \beta)-a(1+2n-\beta)b_n+a b_n^2+2a
x_n-x_n^2}{a(x_n-a)}.
\end{equation}
Differentiating equation (\ref{bn'}) and using the expression for
$x_n'$, we get an equation for $b_n'',\,b_n',\,b_n$ and $x_n$.
Finally, we can eliminate $x_n$ by computing the resultant of this
equation and equation (\ref{bn'}) and we get a nonlinear second
order and second degree equation
\begin{equation}\label{nonlinear diff}
G(b_{n}'',b_n',b_{n},a)=0.
\end{equation}
Now the difficulty is in identifying this equation. Let us scale
the independent variable $b_n(a)=B_n(a/k_1)$ and take $a=k_1 t$.
The equation (\ref{nonlinear diff}) becomes
$G_1(B_n'',B_n',B_n,t)=0,$ where the differentiation is with
respect to $t$. The following result holds.

\begin{theorem}
The equation $G_1(B_n'',B_n',B_n,t)=0$ is reduced to the fifth
Painlev\'e equation \pf\ by the following transformation
\begin{equation}\label{transf}
B_n(t)=\frac{1+n+(\beta-3n-2)y+(1+2n-\beta)y^2+t y'}{2y(y-1)},
\end{equation}
where $y=y(t)$ satisfies \pf\ with parameters   given by
\begin{equation}\label{P5 par}
A=\frac{(\beta-1)^2}{2},\;\;B=-\frac{(n+1)^2}{2},\;\;C=2k_1,\;\;D=0.
\end{equation}
\end{theorem}

\begin{remark}
It is known \cite{GrLSh} that \pf\ with $D=0$ can be reduced to
the third Painlev\'e equation. However, to show that equation
(\ref{nonlinear discrete}) can  be obtained from the B\"acklund
transformation, it is more convenient to use \pf.
\end{remark}

\begin{remark}
The parameters (\ref{P5 par}) are invariant under
$\beta\rightarrow 2-\beta$; compare with the weight on the shifted
lattice (\ref{weight sh}).
\end{remark}

Let $k_1=1$ and $a=t$. We use Theorem \ref{thm: BTr} to get
$y_{n+1}$ and $y_{n-1}$ from $y=y_n.$ In particular,
$$y_{n-1}=1-\frac{4t(y-1)y^2}{(\beta-1)^2(y-1)^2y^2-((n+1)(y-1)-t
y')^2}$$ is a  solution of \pf\ with
$A=(\beta-1)^2/2,\;B=-n^2/2,\;C=2,\;D=0$ and
$$y_{n+1}=1-\frac{4t(y-1)y^2}{(\beta-1)^2(y-1)^2y^2-((n+1)(y-1)+ty')^2}$$
is a solution of \pf\ with
$A=(\beta-1)^2/2,\;B=-(n+2)^2/2,\;C=2,\;D=0$ provided that
$y=y(t)$ is a solution of \pf\ with (\ref{P5 par}). Using
(\ref{transf}) we can express $b_{n\pm 1}$ in terms of $y',\,y$.
For instance,
$$b_{n+1}=\frac{1+n+(\beta-3n-4)y+(3+2n-\beta)y^2-t
y'}{2y(y-1)}.$$ Substituting (\ref{transf}) and the expressions of
$b_{n\pm 1}$ into (\ref{nonlinear discrete}), we indeed see that
they satisfy the equation. Therefore, the system
(\ref{discrete1})--(\ref{discrete2}) is related to the B\"acklund
transformation of \pf.

Finally we use Theorem \ref{P5P3} to obtain the parameters of \pt.
Taking $k_1=1/2$ from (\ref{P5 par}) we get the following values
for different choices of the square roots:
\begin{equation}\label{P3 par1}
\tilde{\alpha}= 2(1+n-\beta),\; \tilde{\beta}=-2(n+\beta),
\end{equation}
\begin{equation}\label{P3 par2}
\tilde{\alpha}= -2(3+n-\beta),\; \tilde{\beta}=2(n+\beta),
\end{equation}
\begin{equation}\label{P3 par3}
\tilde{\alpha}= 2(n-1+\beta),\; \tilde{\beta}=-2(2+n-\beta),
\end{equation}
\begin{equation}\label{P3 par4}
\tilde{\alpha}= -2(1+n+\beta),\; \tilde{\beta}=2(2+n-\beta).
\end{equation}
Taking $k_1=1/2,\;\;a=t/2$ we change variables $t\rightarrow
\sqrt{2t}=z$ in (\ref{nonlinear diff}) and take
$b_n(a)=\tilde{B}_n(z)$. Then the functions
$$\tilde{B}_n(z)=-\frac{z+u(2\beta-2n-1+z u)+z u'}{4u},$$
$$\tilde{B}_n(z)=\frac{1}{4}\left(1+2n-2\beta+z\left(\frac{8(n+1)}{z-u(zu+2\beta-2n-5)-z u'}-\frac{u'+u^2+1}{u}\right)\right)$$
$$\tilde{B}_n(z)=-\frac{z+u(2\beta-2n-1+z u)+z u'}{4u}$$
$$\tilde{B}_n(z)=\frac{1}{4}\left(1+2n-2\beta+z\left(\frac{8(n+1)}{z+u(1+2n+2\beta-zu)-zu'}-\frac{1+u^2+u'}{u}\right)\right)$$
give the transformations to the third Painlev\'e equation for the
function $u=u(z)$ with parameters (\ref{P3 par1})--(\ref{P3 par4})
respectively.

In case $\beta=1$ we get from (\ref{P3 par1})--(\ref{P3 par4})
either
$$\tilde{\alpha}=2n,\;\tilde{\beta}=-2(n+1)$$ or
$$\tilde{\alpha}=-2(n+2),\;\tilde{\beta}=2(n+1).$$

\subsection{Initial conditions}
In this section we study the initial conditions (\ref{initial}),
(\ref{initial sh}), (\ref{initial bi}) for the generalized
Charlier weight (\ref{weight}) on the lattice $\N$, on the shifted
lattice $\N+1-\beta$ and on the bi-lattice respectively. We show
that the initial conditions correspond to classical solutions of
\pf.

First we get that the equation (\ref{nonlinear diff}) with  $n=0$
has solutions satisfying
\begin{equation}\label{Ric b} a^2 b_0'+a
b_0^2-a(1-\beta)b_0-a^2=0,
\end{equation}
where $'=d/da.$  The initial conditions (\ref{initial}),
(\ref{initial sh}), (\ref{initial bi}) satisfy this equation.
Moreover, by using the formulas from \cite{AbramowitzStegun, NIST}, or
computational software,  it can be shown that the general solution
of (\ref{Ric b}) is given by (\ref{initial bi}). We note that  the
Toda system (\ref{Toda}) and discrete system system
(\ref{discrete1})--(\ref{discrete2}) are unchanged under the
transformation $\beta\rightarrow 2-\beta,$
$a_n^2=\hat{a}_n^2,\;\;b_n=\hat{b}_n-1+\beta$.  Also equation
(\ref{Ric b}) is unchanged after this transformation with $n=0$.

Using the transformation (\ref{transf}) with $k_1=1,\;\;a=t$ in
(\ref{Ric b}) we get a first order second degree equation
$$(y-1)(1+y((\beta-1)^2y^2-(4t+(\beta-1)^2)y-1))+2t(y-1)y'-t^2
y'^2=0,$$ for which the solutions satisfy \pf\ with parameters
$A=(\beta-1)^2/2,\;B=-1/2,\;C=2,\;D=0$. Using the transformation
$y\rightarrow 1/y$ in \pf\ we get equation (38.7) in \cite{GrLSh}
for special values of the parameters. Hence, the solutions of \pf\
in terms of classical functions correspond to the initial
conditions for the recurrence coefficients.

\begin{remark}{\rm The recurrence coefficients $a^2_n$ and $b_n$ can
always be written \cite{Chihara} as ratio's of Hankel determinants
containing the moments of the orthogonality measure. However, the
explicit determinant formulas of classical solutions for  all the
Painlev\'e equations are  known.}
\end{remark}

\subsection{Revisiting case $\beta=1$}

It is shown in \cite{LiesGFWalter} that the monic polynomials
orthogonal with respect to a semi-classical variation of the
Charlier weight given by
$$w(k)=\frac{a^k}{(k!)^2},\qquad a>0,\;\;k\in\N,$$ have recurrence
coefficients which satisfy \pf\ with
\begin{equation}\label{par2}A=-B=\frac{n^2}{8},\;\;C=0,\;\;D=-8.\end{equation} This weight can be
obtained from weight (\ref{weight}) by  taking $\beta=1$. However,
if we take the parameter $\beta=1$ in (\ref{P5 par}), we get \pf\
with
\begin{equation}\label{par1}
A=0,\;\;B=-\frac{(n+1)^2}{2},\;\;C=2k_1,\;\;D=0,
\end{equation}
which seems to be a contradiction. In the following we will show how to reduce (\ref{nonlinear
diff}) to \pf\ with parameters (\ref{par2}).

The following transformation exists  between solutions of \pf\
\cite{GrLSh}. Let $y=y(z)$ be a solution of \pf\ with
$A=-B=\alpha/4,\;\;C=0,\;\;D=4\gamma$. Then the function
\begin{equation}\label{transf1}\tilde{y}(t)=\frac{(y+1)^2}{4y}, \;\;t=z^2\end{equation}
is a solution of \pf\ with
$A=\alpha,\;\;B=0,\;\;C=\gamma,\;\;D=0.$ So, if we start with a
solution $y$ of \pf\ with parameters (\ref{par1}), we can first
use the transformation $y\rightarrow 1/y$ to get the parameters
$$A=\frac{(n+1)^2}{2},\;\;B=0,\;\;C=-2k_1,\;\;D=0,$$ then the
transformation (\ref{transf1}) to get solutions with
$$A=\frac{(n+1)^2}{8},\;\;B=-\frac{(n+1)^2}{8},\;\;C=0,\;\;D=-8k_1.$$
Then taking $k_1=1 $ and using Theorem \ref{thm: BTr2} (e.g., a
transformation from the example in case $\varepsilon=1$) we get a
solution of \pf\ with parameters (\ref{par2}). Explicitly this
gives
$$b_n(t)=\frac{n+7ny^2-ny^3-2z y'-y(7n+2z
y')}{8y(y-1)},\;\;t=z^2,$$ where $y=y(z)$ is a solution of \pf\
with (\ref{par2}).  Equation (\ref{Ric b}) is also transformed to
a first order second degree equation for which the solutions
satisfy \pf. Hence, we get that the solutions of \pf\ in terms of
classical special functions generate the initial conditions for the
system (\ref{discrete1})--(\ref{discrete2}).

\section{Discussion}

The recurrence coefficients for semi-classical weights are  often
related to the solutions of the  Painlev\'e type equations (e.g.,
\cite{GFWalterLaguerre, GFWalterMeixner} and see also  overview in
\cite{LiesGFWalter}). The known results and our recent findings of
the relations between the recurrence coefficients of the
(generalized) Charlier and Meixner polynomials   and the fifth
Painlev\'e equation (\ref{P5}) with parameters $(A,B,C,D)$ can be
summarized as follows. We assume below that $
k\in\mathbb{N}=\{0,1,\ldots\}$.

\begin{enumerate}

\item \begin{enumerate}

\item The weight $a^k/k!,\;a>0, $ is a classical Charlier weight and the
recurrence coefficients can be computed explicitly.

\item  It is shown in \cite{LiesGFWalter} that the
recurrence coefficients for the weight $a^k/(k!)^2, \;\;a>0,$ are
related to the solutions of the  fifth Painlev\'e equation \pf\
with parameters $(n^2/8,-n^2/8,0,-8)$. Moreover, as shown in this
paper, by using certain transformations, they are also related to
solutions of  \pf\ with parameters $(0,-(n+1)^2/2,2,0)$ (this is
the special case of 1(c)).

\item As shown in this paper, the
recurrence coefficients for the weight $a^k/((\beta)_k
k!),\;\;a>0,\;\beta\neq 1,$ are related to solutions of  \pf\ with
$((\beta-1)^2/2,-(n+1)^2/2,2,0)$.

\end{enumerate}
\item \begin{enumerate}

\item The weight $(\beta)_k c^k/k!,\;\;\beta>0,\;0<c<1,$ is the classical Meixner
weight and the recurrence coefficients are known explicitly.

\item   It is shown in   \cite{LiesGFWalter} that the
recurrence coefficients for  the weight $(\beta)_k
c^k/(k!)^2,\;\;\beta>0,\;c>0,$ are related to solutions of  \pf\
with $((\beta-1)^2/2,-(\beta+n)^2/2,2n,-2)$.

\item It is shown in \cite{GFWalterMeixner} that
the recurrence coefficients for the weight $(\gamma)_k
c^k/(k!(\beta)_k),\;\;c,\beta,\gamma>0,$ are related to solutions
of  \pf\ with
$((\gamma-1)^2/2,-(\gamma-\beta+n)^2/2,k_1(\beta+n),-k_1^2/2))$,
$k_1\neq 0$. The case $\gamma=1$ corresponds to  classical
Charlier weight  on the shifted lattice $\mathbb{N}+1-\beta$. The
case $\beta=\gamma$ gives the  Charlier weight  and the case
$\beta=1$ is given in 2(b).

\end{enumerate}

\end{enumerate}

Therefore, we have essentially two different cases for \pf: $D=0$
for the recurrence coefficients of the  generalizations of the
Charlier weight and $D\neq 0$ for the recurrence coefficients of
the generalizations of the Meixner weight. Note that the case
$D=0$ can be reduced to the third Painlev\'e equation
\cite{GrLSh}.

It is  interesting to study weights on shifted lattices. For
instance, there is no shifted lattice for the weight
$a^k/(k!)^2 $  since both lattices   coincide. Since  the lattices
are determined by the zeros of the weight,  one could have more
than one lattice if one takes $a^k/((\beta)_k(\gamma)_k(\delta)_k
k!)$, which is zero for $k=-1$, $k=-\beta$, $k=-\gamma$ and
$k=-\delta$.

It is an interesting open problem  to study recurrence
coefficients for more complicated weights  like $a^k/(k!)^3$ or
$(\gamma)_k(\beta)_k c^k/(k!)^2$ and others. The Toda system will
be the same but the discrete system will be much more complicated
(maybe higher order discrete Painlev\'e equations will  be
needed). Will this lead to the Painlev\'e equations  (with more
general parameters $A,B,C,D$ depending on the parameters in the
weight) or some higher order equations in some hierarchy of the
Painlev\'e equations? The direct calculations will be cumbersome
and other techniques like using the Lax pairs or others will be
needed to tackle such problems.

\section*{Acknowledgements} GF is
 partially supported by Polish MNiSzW Grant N N201 397937.  WVA is supported by Belgian Interuniversity Attraction Pole P6/02, FWO
grant G.0427.09 and K.U.Leuven Research Grant OT/08/033.

\begin{verbatim}
Galina Filipuk
Faculty of Mathematics, Informatics and Mechanics
University of Warsaw
Banacha 2 Warsaw 02-097, Poland
filipuk@mimuw.edu.pl

Walter Van Assche
Department of Mathematics
Katholieke Universiteit Leuven
Celestijnenlaan 200B box 2400
BE-3001 Leuven, Belgium
Walter.VanAssche@wis.kuleuven.be
\end{verbatim}


\begin{thebibliography}{99}

\bibitem{AbramowitzStegun} M. Abramowitz, I. Stegun,
{\it Handbook of Mathematical Functions}, Dover Publications, New
York, 1965.

\bibitem{Adler} V.E.  Adler, {\it Nonlinear chains and Painlev\'e equations}, Physica D {\bf 73} (1994), 335--351.

\bibitem{LiesGFWalter} L. Boelen, G. Filipuk and W. Van Assche, {\it Recurrence coefficients of generalized Meixner
polynomials and Painlev\'e equations}, J. Phys. A: Math. Theor.
{\bf 44} (2011), 035202 (19 pp).

\bibitem{Chihara} T. S. Chihara, {\it An Introduction to Orthogonal Polynomials}, Gordon and Breach, New York,
1978.

\bibitem{GFWalterLaguerre} G. Filipuk, W. Van Assche, L. Zhang {\it The recurrence coefficients of semi-classical Laguerre polynomials and the fourth Painlev\'e equation}, \texttt{arXiv:1105.5229v1 [math.CA]}

\bibitem{GFWalterMeixner} G. Filipuk and W. Van Assche, {\it   Recurrence coefficients of a new generalization of
the Meixner polynomials}, preprint \texttt{arXiv:1104.3773v1 [math.CA]}.

\bibitem{GrLSh} V.I. Gromak, I. Laine and S. Shimomura,
{\it Painlev\'e  Differential Equations in the Complex Plane},
Vol.28, Studies in Mathematics, de Gruyter, Berlin, NewYork, 2002.

\bibitem{Noumi} M. Noumi, {\it  Painlev\'e equations through
symmetry},   Translations of Mathematical Monographs, {\bf 223},
American Mathematical Society, Providence, RI, 2004.

\bibitem{NIST} F.W.J. Olver, D.W. Lozier, R.F. Boisvert, C.W. Clark,
{\it NIST Handbook of Mathematical Functions}, NIST and Cambridge University Press, 2010.

\bibitem{bilattice} C. Smet and W. Van Assche, {\it Orthogonal polynomials on a
bi-lattice}, Constr. Approx. (to appear); \texttt{arXiv:1101.1817v1 [math.CA]}.

\bibitem{Tsegelnik} V.V.  Tsegelnik, {\it The Painlev\'e type equations:
analytical properties of solutions and their applications}, Habilitation thesis, Minsk, 2001 (in Russian).

\end{thebibliography}
\end{document}